\documentclass{article}

\usepackage{
 amsmath,
 amsxtra,
 amsthm,
 amssymb,
 etex,
 mathrsfs,
 }
\usepackage[all]{xy}

\theoremstyle{plain}
\newtheorem{thm}{Theorem}[subsection]
\newtheorem{cor}[thm]{Corollary}
\newtheorem{lem}[thm]{Lemma}
\newtheorem{prop}[thm]{Proposition}
\newtheorem*{thm2}{Theorem}

\theoremstyle{definition}
\newtheorem{defn}[thm]{Definition}
\newtheorem*{remark}{Remark}

\newcommand{\bd}{\begin{defn}}
\newcommand{\ed}{\end{defn}}
\newcommand{\bl}{\begin{lem}}
\newcommand{\el}{\end{lem}}
\newcommand{\bp}{\begin{prop}}
\newcommand{\ep}{\end{prop}}
\newcommand{\bt}{\begin{thm}}
\newcommand{\et}{\end{thm}}
\newcommand{\bc}{\begin{cor}}
\newcommand{\ec}{\end{cor}}
\newcommand{\br}{\begin{remark}}
\newcommand{\er}{\end{remark}}
\newcommand{\bdi}{\begin{diagram}}
\newcommand{\edi}{\end{diagram}}
\newcommand{\beq}{\begin{equation}}
\newcommand{\eeq}{\end{equation}}
\newcommand{\ba}{\begin{array}}
\newcommand{\ea}{\end{array}}
\newcommand{\bpf}{\begin{proof}}
\newcommand{\epf}{\end{proof}}

\newcommand{\Z}{\mathbb{Z}}
\newcommand{\Q}{\mathbb{Q}}
\newcommand{\Zp}{\mathbb{Z}_{p}}
\newcommand{\Qp}{\mathbb{Q}_{p}}

\newcommand{\al}{\alpha}

\newcommand{\Ga}{\Gamma}
\newcommand{\ga}{\gamma}

\newcommand{\la}{\lambda}
\newcommand{\si}{\sigma}
\newcommand{\Hi}{H_{\mathrm{Iw}}}

\DeclareMathOperator{\Sel}{Sel} \DeclareMathOperator{\Gal}{Gal}
 
\DeclareMathOperator{\corank}{corank}
\DeclareMathOperator{\Ext}{Ext}

\newcommand{\cyc}{\mathrm{cyc}}
\newcommand{\M}{\mathfrak{M}}

\newcommand{\ot}{\otimes}
\newcommand{\ilim}{\displaystyle \mathop{\varinjlim}\limits}
\newcommand{\plim}{\displaystyle \mathop{\varprojlim}\limits}
\newcommand{\im}{\mathrm{im}\,}

\newcommand{\lra}{\longrightarrow}

\newcommand{\ps}[1]{[[ #1 ]]}

\begin{document}
\title{On the signed Selmer groups of congruent elliptic curves with semistable reduction at all primes above $p$}
 \author{ Suman Ahmed\footnote{School of Mathematics and Statistics,
Central China Normal University, Wuhan, 430079, P.R.China.
 E-mail: \texttt{npur.suman@gmail.com}}  \quad
  Meng Fai Lim\footnote{School of Mathematics and Statistics $\&$ Hubei Key Laboratory of Mathematical Sciences,
Central China Normal University, Wuhan, 430079, P.R.China.
 E-mail: \texttt{limmf@mail.ccnu.edu.cn}} }
\date{}
\maketitle

\begin{abstract} \footnotesize
\noindent Let $p$ be an odd prime. We attach appropriate signed Selmer groups to an elliptic curve $E$, where $E$ is assumed to have semistable reduction at all primes above $p$. We then compare the Iwasawa $\lambda$-invariants of these signed Selmer groups for two congruent elliptic curves over the cyclotomic $\Zp$-extension in the spirit of Greenberg-Vatsal and B. D. Kim. As an application of our comparsion formula, we show that if the $p$-parity conjecture is true for one of the congruent elliptic curves, then it is also true for the other elliptic curve. In the midst of proving this latter result, we also generalize an observation of Hatley on the parity of the signed Selmer groups.


\end{abstract}

\section{Introduction}

Let $\Q^{\cyc}$ be the cyclotomic $\Zp$-extension of $\Q$, where $p$ is a fixed odd prime. In \cite{GV}, Greenberg and Vatsal showed that if two congruent elliptic curves have good ordinary reduction at $p$, then the Iwasawa $\la$-invariants of their classical $p$-primary Selmer groups are related in an explicit manner.
When the elliptic curve has supersingular reduction, the classical $p$-primary Selmer groups are not torsion over the Iwasawa algebra, and it was due to the insight of Kobayashi \cite{Kob} that one has to study the so-called signed Selmer groups. Building off on this, B. D. Kim \cite{Kim09} established a relation for the Iwasawa $\la$-invariants of the signed Selmer groups of two congruent elliptic curves.


The goal of this paper is to extend the results of Greenberg-Vatsal and B. D. Kim to the situation, where the elliptic curves have no additive reduction at all primes above $p$. In other words, our elliptic curves, which are now defined over a finite extension of $\Q$, can have either good ordinary/supersingular
reduction or split/non-split multiplicative reduction at each prime above $p$. We shall briefly describe our results here.

As a start, we introduce the cast. Suppose that $E_1$ and $E_2$ are two elliptic curves defined over a number field $L'$ such that $E_1[p]\cong E_2[p]$ as $\Gal(\overline{L'}/L)$-modules, where $L$ is a finite extension of $L'$.  As our results will be proved over $L$ and $L(\mu_p)$, we shall always write $F$ to denote either one.

We follow up by introducing the following standing assumptions on our cast.

\begin{enumerate}
\item[$\mathbf{(NA)}$] The elliptic curves $E_1$ and $E_2$ have no additive reduction at all primes of $L'$ above $p$.

\item[$\mathbf{(SS)}$]  For each prime $u$ of $L'$ above $p$ which is a good supersingular reduction prime of $E_1$ and $E_2$, the prime $u$ is unramified in $L/L'$ with $L'_u=\Qp$ and $a_u = 1 + p - |\tilde{E}_{i,u}(\mathbb{F}_p)| = 0$, where $\tilde{E}_{i,u}$ is the reduction of $E_i$ at $u$ ($i=1,2)$.

\item[$\mathbf{(RED)}$] The elliptic curves $E_1$ and $E_2$ have the same reduction type at each prime of $L'$ above $p$. Furthermore, in the event that $E_1$ and $E_2$ have non-split multiplicative reduction at a prime $u$ of $L'$ above $p$, suppose further that $E_1$ and $E_2$ have the same reduction type at all primes of $F$ above $u$, where we recall that $F$ denotes either $L$ or $L(\mu_p)$.

\item[$\mathbf{(RU)}$] For each prime $u$ of $L'$ above $p$ which is either a prime of good ordinary reduction or non-split multiplicative reduction of $E_1$ (and hence $E_2$ in view of $\mathbf{(RED)}$), assume that $L'_u$ does not contain $\mu_p$.
\end{enumerate}

Write $F^\cyc$ for the cyclotomic $\Zp$-extension of $F$. Following \cite{Kim07,Kim09,KimPM,Kim14,KO,Kob}, we define a mixed signed Selmer group $\Sel^{\overrightarrow{s}}(E_i/F^{\cyc})$ of $E_i$ over $F^{\cyc}$ (see Subsection \ref{subsec signed Selmer} for the precise definition). Each of these Selmer groups has a natural $\Zp\ps{\Ga}$-module structure, where $\Ga= \Gal(F^\cyc/F)$.  Denote by $\la^{\overrightarrow{s}}_{E_i}$ the Iwasawa $\la$-invariant of $\Sel^{\overrightarrow{s}}(E_i/F^{\cyc})$ for $i=1,2$.

It is well-known from the structure theory of $\Zp\ps{\Ga}$-modules (cf. \cite[Section 5.3]{NSW}) that a finitely generated $\Zp\ps{\Ga}$-module $M$ is finitely generated over $\Zp$ if and only if $M$ is a torsion $\Zp\ps{\Ga}$-module with trivial Iwasawa $\mu$-invariant. Furthermore, in the event of such, the Iwasawa $\la$-invariant of $M$ coincides with the $\Zp$-rank of $M$.
In particular, if
$\Sel^{\overrightarrow{s}}(E_i/F^{\cyc})$ is cofinitely generated over $\Zp$, we have $\la^{\overrightarrow{s}}_{E_i} = \corank_{\Zp}\Big(\Sel^{\overrightarrow{s}}(E_i/F^{\cyc})\Big)$. The Iwasawa invariants of these signed Selmer groups are intimately related to the growth of the Mordell-Weil groups and Tate-Shafarevich groups of the elliptic curve in the cyclotomic $\Zp$-extension (for instance, see \cite[Proposition 4.4]{KO}, \cite[Corollary 10.2 and Theorem 10.9]{Kob} and \cite[Theorem 2]{LLSha}).

Our first main theorem is as follows.

\begin{thm2}[Theorem \ref{congruent theorem}]
Retain the above settings. In particular, suppose that $\mathbf{(NA)}$, $\mathbf{(SS)}$, $\mathbf{(RED)}$ and $\mathbf{(RU)}$ are valid.
Then $\Sel^{\overrightarrow{s}}(E_1/F^{\cyc})$ is cofinitely generated over $\Zp$ if and only if
$\Sel^{\overrightarrow{s}}(E_2/F^{\cyc})$ is cofinitely generated over $\Zp$. Moreover, if this is so, we have the following equality
\[ \la^{\overrightarrow{s}}_{E_1} + \sum_{w\in S'(F^{\cyc})}\corank_{\Zp}\Big(H^1(F^{\cyc}_w,E_1(p))\Big)= \la^{\overrightarrow{s}}_{E_2} + \sum_{w\in S'(F^{\cyc})}\corank_{\Zp}\Big(H^1(F^{\cyc}_w,E_2(p))\Big). \]
\end{thm2}

We say a little bit on the standing assumptions in the theorem. Hypothesis $\mathbf{(SS)}$ is necessary for the definition of the signed Selmer groups, as this allows one to define the plus/minus norm subgroups of Kobayashi. Since the proof of the theorem requires comparing the local cohomology groups for primes above $p$, hypotheses $\mathbf{(RED)}$ and $\mathbf{(RU)}$ are imposed to allow us to make such meaningful comparison on the local cohomology groups. The hypotheses $\mathbf{(RED)}$ and $\mathbf{(RU)}$ are also innate in \cite{AAS, GV, Hat, Kim09, She}. In fact, in these works, their congruent elliptic curves are always assumed to have only good ordinary reduction or only good supersingular reduction above $p$, and so $\mathbf{(RED)}$ automatically holds. In \cite{AAS, GV, She}, as their elliptic curves are defined over $\Q$, $\mathbf{(RU)}$ is automatically satisfied. It is
not too difficult to obtain examples of elliptic curves with mixed reduction types at primes above $p$ by arguments similar to those in \cite[Proposition 5.4]{G99} or \cite[Lemma 8.19]{Maz}, although we have to confess that we do not know how to compute the Iwasawa invariants on any of these examples at this point of the writing.

Note that even in the case where the elliptic curves have ordinary reduction at all primes above $p$, our result is an improvement of the results of Greenberg-Vatsal \cite{GV} and Ahmed-Aribam-Shekhar \cite{AAS, She} for we do not require the assumption that $E(F)[p]=0$.
We should mention that the results of Greenberg-Vatsal and Kim have also been established for modular forms of higher weight and even more general Galois
representations (for instances, see \cite{EPW, Hac, HatLei, Pon}). However in these prior works, they always work with coherent reduction types above $p$. In other words, their Galois representations are assumed to have either ordinary reduction at all primes above $p$ or non-ordinary reduction at all primes above $p$.

As an application of Theorem \ref{congruent theorem}, we investigate the $p$-parity conjecture following \cite{AAS, Hat, She}.  For an elliptic curve $E$ over $F$, the $p$-parity conjecture predicts that $w(E/F) = (-1)^{s_p(E)}$ (see \cite{Do} and references loc. cit), where we write $w(E/F)$ for its global root number and $s_p(E)$ for the $\Zp$-corank of $\Sel(E/F)$.

We require another hypothesis which is required in handling signed Selmer groups, when one of the signs is a $+$ sign.

\medskip
$\mathbf{(S+)}$ For every $v\in S_p^{ss}(L)$, assume further that $[L_v:\Qp]$ is not divisible by 4.

\smallskip
 Under this hypothesis, one can show that the $\la$-invariants of the signed Selmer groups have the same parity as the $\Zp$-corank of the classical Selmer group $\Sel(E_i/F)$ (see Proposition \ref{selmer parity}). Combining this with Theorem \ref{congruent theorem} and with a bit more work, we obtain the following.

\begin{thm2}[Theorem \ref{p-parity theorem}]
Retain the settings as above. Suppose further that at least one of the following statements is also valid.
\begin{enumerate}
\item[$(a)$] $\Sel^{\overrightarrow{-}}(E_1/F^{\cyc})$ is cofinitely generated over $\Zp$.
\item[$(b)$] $\Sel^{\overrightarrow{s}}(E_1/F^{\cyc})$ is cofinitely generated over $\Zp$ for some $\overrightarrow{s}\neq \overrightarrow{-}$ and $(\mathbf{S+})$ is valid.
    \end{enumerate}
 Then one has
\[ \frac{w(E_1/F)}{w(E_2/F)} = (-1)^{s_p(E_1) - s_p(E_2)}. \]
In particular, the $p$-parity conjecture holds for $E_1$ over $F$ if and only if it holds for $E_2$ over $F$.
\end{thm2}

When the elliptic curves either have good ordinary reduction at all primes above $p$ or good supersingular reduction at all primes above $p$, the above result was proved under the assumption of statement (a) by the first named author of this paper with Aribam and Shekhar in \cite[Corollary 5.6]{AAS}. Our result is an improvement of this, where we allow mixed reduction types above $p$ and mixed signs in the definition of our signed Selmer groups. We should also mention that the proof of  \cite[Corollary 5.6]{AAS} contains a slight gap which we have addressed in this paper (see Remark after Proposition \ref{selmer parity}).

We now give an outline of our paper. In Section \ref{local calculations}, we collect various results concerning the arithmetic of an elliptic curve over a local field which will be needed for the subsequent discussion of the paper. Section \ref{Selmer} is where we introduce the signed Selmer group of an elliptic curve and give various equivalent descriptions of them. The point of having these different descriptions is that in proving Theorem \ref{congruent theorem}, it is more natural to work with the description of the signed Selmer groups as given in Proposition \ref{selmer alternative}. On the other hand, while proving Theorem \ref{p-parity theorem}, one needs to work with the so-called strict signed Selmer groups which coincides with the signed Selmer groups on the level of $F$ and $F^{\cyc}$. We do however remark that the strict signed Selmer groups and the signed Selmer groups need not coincide on the intermediate subextensions of $F^{\cyc}/F$ (see the proof of Proposition \ref{strict selmer=selmer}). Finally, Section \ref{main results} is where all our main results will be established. In the midst of proving our second main theorem, we also generalize an observation of Hatley \cite[Corollary 4.2]{Hat} on the parity of the signed Selmer groups (see Corollary \ref{selmer parity}). Although this latter result is not required for the proof of our second main theorem, we have thought that it is interesting enough to be noted down (also see Remark after Corollary \ref{selmer parity}).

\subsection*{Acknowledgement}
We would like to thank Antonio Lei and Ramdorai Sujatha for making us aware the validity of Proposition \ref{torsion surjective H2} and for sharing their work \cite{LeiSu} with us. We would also like to thank the referee for various helpful comments which helped us in improving the exposition of the paper.
The research of this article took place when S. Ahmed was a postdoctoral fellow at Central China Normal University, and he would like to acknowledge the hospitality
and conducive working conditions provided by the said institute. M. F. Lim's research is supported by the
National Natural Science Foundation of China under Grant No. 11550110172 and Grant No. 11771164.

\section{Local consideration} \label{local calculations}

In this section, we collect certain results on elliptic curves defined over a local field.

\subsection{Ordinary reduction at $p$} \label{ordinary subsection}

Let $K$ be a finite extension of $\Q_p$ and $E$ an elliptic curve defined over $K$. In this subsection, our elliptic curve $E$ is always assumed to have either good ordinary reduction or multiplicative reduction. Then from \cite[P. 150]{CG}, we have the following short exact sequence
\[ 0\lra C \lra E(p)\lra D\lra 0\]
of $\Gal(\bar{K}/K)$-modules, where $C$ and $D$ are cofree $\Zp$-modules of corank one, and are characterized by the facts that $C$ is divisible and that $D$ is the maximal quotient of $E(p)$ by a divisible subgroup such that $\Gal(\bar{K}/K^{ur})$ acts on $D$ via a finite quotient. Here $K^{ur}$ is the maximal unramified extension of $K$. In fact, $D$ can be explicitly described as follows (see \cite{CG} or \cite[Section 3]{G99})
\[D = \begin{cases}  \widetilde E(p),& \mbox{if $E$ has good ordinary reduction}, \\
      \Qp/\Zp, & \mbox{if $E$ has split multiplicative reduction,} \\
      \Qp/\Zp\ot\phi, & \mbox{if $E$ has nonsplit multiplicative reduction,}\end{cases}
 \]
where $\widetilde E$ is the reduction of $E$ and $\phi$ is a unramified character of $\Gal(\bar{K}/K)$. Note that in all these three cases, $D$ is a trivial $\Gal(\bar{K}/K^{ur})$-module which in turn enables us to view $D$ as a $\Gal(K^{\cyc,ur}/K^{\cyc})$-module. We now break our discussion into three cases.






Case 1: Suppose that $E$ has good ordinary reduction. Since $D(K^{\cyc}) = \widetilde{E}(K^{\cyc})$ is finite, the group $\Gal(K^{\cyc,ur}/K^{\cyc})$ does not act trivially on $D$. Consequently, $D_{\Gal(K^{\cyc,ur}/K^{\cyc})}$ is a proper quotient of $D$. But since $D$ is divisible of $\Zp$-corank one, this forces $D_{\Gal(K^{\cyc,ur}/K^{\cyc})}=0$, or equivalently, $H^1\big(\Gal(K^{\cyc,ur}/K^{\cyc}), D\big) = 0$. On the other hand, we clearly have
$H^2\big(\Gal(K^{\cyc,ur}/K^{\cyc}), D\big) = 0$ as $\Gal(K^{\cyc,ur}/K^{\cyc})\cong \hat{\Z}$ has $p$-cohomological dimension one. It then follows from the inflation-restriction sequence that $H^1(K^{\cyc},D)\cong H^1(K^{\cyc, ur}, D)$.

Case 2: Suppose that $E$ has non-split multiplicative reduction. Then $D=\Qp/\Zp\ot\phi$ for a unramifed character of $\Gal(\bar{K}/K)$ which factors through a quadratic extension of $K$. Since $p$ is odd, the multiplicative reduction remains non-split. Hence $D(K^{\cyc})$ is finite. A similar argument as above then yields an isomorphism $H^1(K^{\cyc},D)\cong H^1(K^{\cyc, ur}, D)$.

Case 3: Suppose that $E$ has split multiplicative reduction. Then we have $D\cong \Qp/\Zp$ as $\Gal(\bar{K}/K^{\cyc})$-modules.

We record the above analysis into the following lemma which will be required in the subsequent of the paper.

\bl \label{ordinary}
Let $E$ be an elliptic curve defined over a finite extension $K$ of $\Qp$. Then the following statements hold.
\begin{enumerate}
\item[$(a)$] If $E$ has good ordinary reduction or non-split multiplicative reduction, then $D(K^{\cyc})$ is finite and there is an isomorphism $H^1(K^{\cyc},D)\cong H^1(K^{\cyc, ur}, D)$.
\item[$(b)$] If $E$ has split multiplicative reduction, then $D\cong \Qp/\Zp$ as $\Gal(\bar{K}/K^{\cyc})$-modules.
\end{enumerate}
\el

We end the subsection with the following lemma concerning a uniqueness property of $C[p]$ which will be required for the proof of Theorem \ref{congruent theorem}.

\bl \label{D[p]}
Suppose that the local field $K$ does not contain $\mu_p$. Then $C[p]$ is the unique $\Gal(\bar{K}/K)$-submodule of $E[p]$ which is isomorphic to $\mu_p$ as $\Gal(\bar{K}/K^{ur})$-modules.
\el

\bpf
Let $U$ be a $\Gal(\bar{K}/K)$-submodule of $E[p]$ which is isomorphic to $\mu_p$ as a $\Gal(\bar{K}/K^{ur})$-module. Now since $U/(U\cap C[p])\cong (U+C[p])/C[p]$ is contained in $E[p]/C[p]\cong D[p]$, it has trivial $\Gal(\bar{K}/K^{ur})$-action. On the other hand, as $K$ does not contain $\mu_p$, $\Gal(\bar{K}/K^{ur})$ does not act trivially on $U$. Therefore, $U/(U\cap C[p])$ must be a proper quotient of $U$. But as $U$ is an one dimensional $\mathbb{F}_p$-vector space, we must have $U/(U\cap C[p])=0$, or equivalently, $U\subseteq C[p]$. Now since both $U$ and $C[p]$ are one dimensional $\mathbb{F}_p$-vector spaces, the inclusion is an equality as required.
\epf

\subsection{Good supersingular reduction at $p$}

In this subsection, $E$ will denote an elliptic curve defined over $\Qp$ with $a_p=0$. Denote by $\hat{E}$ the formal group of $E$. Let $k$ be an unramified extension of $\Qp$ of degree $d$. For $n\geq -1$, set $K_n = k(\mu_{p^{n+1}})$. In particular, $K_{-1} = k$. Note that the extension $K_{\infty}/k$ is totally ramified and $\Gal(K_0/k)\cong \Z/(p-1)\Z$.

We shall write $\hat{E}(K_n) = \hat{E}(\mathfrak{m}_n)$, where $\mathfrak{m}_n$ is the maximal ideal of the ring of integers of $K_n$. Finally, we write $\hat{E}(K_{\infty}) =\cup_n \hat{E}(K_n)$.

\bl \label{supersingular points}
 Retain settings as above. Then $\hat{E}(K_{n})$ and $E(K_n)$ have no $p$-torsion for every $n\geq -1$.
 \el

\bpf
The first assertion follows from \cite[Proposition 3.1]{KO}. Now consider the following short exact sequence
\[ 0\lra \hat{E}(K_n)\lra E(K_n) \lra \widetilde{E}(\kappa_n)\lra 0,\]
where $\kappa_n$ is the residue field of $K_n$. In view of the first assertion
and that $\widetilde{E}(\kappa_n)$ has no $p$-torsion by the supersingular assumption, it follows from the above short exact sequence that $E(K_n)$ have no $p$-torsion as required. \epf

Following \cite{Kim07,Kim09,KimPM,Kim14,KO,Kob}, we define
the following groups
\[\hat{E}^+(K_{n}) = \{ P\in \hat{E}(K_{n})~:~\mathrm{tr}_{n/m+1}(P)\in \hat{E}(K_{m}), 2\mid m, -1\leq m \leq n-1\}, \]
\[\hat{E}^-(K_{n}) = \{ P\in \hat{E}(K_{n})~:~\mathrm{tr}_{n/m+1}(P)\in \hat{E}(K_{m}), 2\nmid m, -1\leq m \leq n-1\}, \]
where $\mathrm{tr}_{n/m+1}: E(K_{n}) \lra E(K_{m+1})$ denotes the trace map.
 We shall then write $\hat{E}^{\pm}(K_{\infty}) =\cup_n \hat{E}^{\pm}(K_n)$.

By fixing a topological generator $\ga$ of $\Gal(K_{\infty}/K_0)$, we identify $\Zp\ps{\Gal(K_{\infty}/k)}$ with the formal power series ring $\big(\Zp[\Gal(K_0/k)]\big)\ps{X}$. For a character $\eta$ of $\Gal(K_0/k)$ and a $\big(\Zp[\Gal(K_0/k)]\big)\ps{X}$-module $M$, let $M^{\eta}$ denotes its $\eta$-eigenspace, which is regarded as a $\Zp\ps{X}$-module. We also write
\[ \delta =\begin{cases} 0,  & \mbox{if $d\neq 0$ (mod 4) or $\eta\neq 1$}, \\
2, & \mbox{otherwise}.\end{cases}\]

With this notation, we can now state the following.

\bp \label{Kita-O results}
Retain settings as above. Write $\Ga_n=\Gal(K_{\infty}/K_n)$. The following statements are then valid.
\begin{enumerate}
\item[$(a)$] $\big(\hat{E}^{\pm}(K_{\infty})^{\eta}\ot\Qp/\Zp\big)^{\Ga_n}$ are cofree $\Zp$-modules for all $n$, and
    \[ \corank_{\Zp}\left(\big(\hat{E}^{s}(K_{\infty})^{\eta}\ot\Qp/\Zp\big)^{\Ga_n}\right)
    =\begin{cases} dp^n+\delta,  & \mbox{if $s=+$}, \\
dp^n, & \mbox{if $s=-$}.\end{cases}\]
\item[$(b)$] One has \[\left(\hat{E}^{s}(K_{\infty})^{\eta}\ot\Qp/\Zp\right)^{\vee}\cong
    \begin{cases} \Zp\ps{X}^{\oplus d}\oplus \Zp^{\oplus \delta},  & \mbox{if $s=+$}, \\
 \Zp\ps{X}^{\oplus d}, & \mbox{if $s=-$}. \end{cases}\]
\end{enumerate}
\ep

\bpf
Statement (a) is \cite[Corollary 3.25]{KO}, and statement (b) follows from \cite[Theorem 3.34]{KO} (also see \cite[Theorems 2.7 and 2.8]{Kim14}).
\epf

Write $\mathbb{H}^{\pm}_{\infty}= \hat{E}^{\pm}(K_{\infty})\ot\Qp/\Zp$ and $\mathbb{H}^{\pm}_n= \left(\hat{E}^{\pm}(K_{\infty})\ot\Qp/\Zp\right)^{\Ga_n}$.
By the Hochshild-Serre spectral sequence and Lemma \ref{supersingular points}, we have an isomorphism
\[ H^1(K_n,E(p))\cong H^1(K_{\infty},E(p))^{\Ga_n}.\]
Via this isomorphism, we may view $\mathbb{H}^{\pm}_n$ as a subgroup of $H^1(K_n,E(p))$.
Let $M_n^{\pm}$ be the exact annihilator of $\mathbb{H}^{\pm}_n$ with respect to the local Tate pairing
\[ H^1(K_n, E(p))\times H^1(K_n, T_pE)\lra \Qp/\Zp. \]
In other words, we have short exact sequences
\[ 0\lra \mathbb{H}^{\pm}_n \lra H^1(K_n, E(p)) \lra \big(M_n^{\pm}\big)^{\vee}\lra 0, \]
\[ 0\lra M_n^{\pm} \lra H^1(K_n, T_pE) \lra \big(\mathbb{H}^{\pm}_n\big)^{\vee}\lra 0. \]
Since each $\mathbb{H}^{\pm}_n$ is divisible by Proposition \ref{Kita-O results}(a), it follows that there are two short exact sequences
\[ 0\lra \mathbb{H}^{\pm}_n[p^j] \lra H^1(K_n, E(p))[p^j] \lra \big(M_n^{\pm}/p^j\big)^{\vee}\lra 0, \]
\[ 0\lra M_n^{\pm}/p^j \lra H^1(K_n, T_pE)/p^j \lra \big(\mathbb{H}^{\pm}_n[p^j]\big)^{\vee}\lra 0. \]

By appealing to Lemma \ref{supersingular points} again, one has $H^1(K_n, E(p))[p^j]\cong H^1(K_n, E[p^j])$ and $H^1(K_n, T_pE)/p^j \cong H^1(K_n, E[p^j])$. Combining this with the above short exact sequences, we see that $M_n^{\pm}/p^j$ is the exact annihilator of $\mathbb{H}^{\pm}_n[p^j]$ with respect to the local Tate pairing
\[ H^1(K_n, E[p^j])\times H^1(K_n, E[p^j])\lra \Z/p^j. \]

In fact, one even has the following.

\bp[Kim] \label{Kim duality}
One always has $\mathbb{H}^{-}_n[p^j] = M_n^{-}/p^j$. If $d$ is not divisible by $4$, then we also have $\mathbb{H}^{+}_n[p^j] = M_n^{+}/p^j$.
\ep

\bpf
The first equality is established by Kim in \cite[Proposition 3.15]{Kim07}. One can check that the same proof carries over for the second equality under the assumption that $d$ is not divisible by $4$ (also see \cite[Theorem 2.9]{Kim14})
\epf

\subsection{Local parity comparison outside $p$} \label{local parity}

In this subsection, $K$ will denote a finite extension of $\Q_l$, where $l\neq p$. Let $E$ be an elliptic curve defined over $K$. Writing $K^{\cyc}$ for the cyclotomic $\Zp$-extension of $K$, it is well-known that $H^1(K^{\cyc}, E(p))$ is cofinitely generated over $\Zp$ (cf. \cite[Proposition 2]{G89}), whose $\Zp$-corank is in turn denoted by $\si_E$. We also write $w(E/K)$ for the local root number of $E$ over $K$ (see \cite[Section 3.4]{Do} for its definition). The following result of Ahmed-Aribam-Shekhar will be required in the subsequent section.

\bp[Ahmed-Aribam-Shekhar] \label{local parity compare} Suppose that $E_1$ and $E_2$ are two elliptic curves defined over $K$ with $E_1[p]\cong E_2[p]$ as $\Gal(\bar{K}/K)$-modules. Then we have
\[ \frac{w(E_1/K)}{w(E_2/K)}=(-1)^{\si_{E_1}-\si_{E_2}}.\]
\ep

\bpf
See \cite[Theorem 5.7]{AAS}.
\epf

\section{Signed Selmer groups over cyclotomic $\Zp$-extension} \label{Selmer}

We now turn to the global situation, where we begin by fixing some notation and standing assumptions that will be adhered throughout
the section. Let $L'$ be a number field and $E$ an elliptic curve defined over $L'$. Fix a finite extension $L$ of $L'$. Let $S$ be a finite set of primes of $L'$ which contains all the primes above $p$, all the ramified primes of $L/L'$, the bad reduction primes of $E$ and the archimedean primes. Denote by $S_p$ the set of primes of $L'$ above $p$. Write $S_p^{ss}$ for the set of primes of $L'$ above $p$ at which $E$ has good supersingular reduction. We also write $S' = S-S_p$ and $S_p^{o} = S_p-S_p^{ss}$. The following assumptions will be in full force throughout this section.

\medskip
$\mathbf{(NA)}$ The elliptic curve $E$ has no additive reduction at all primes of $L'$ above $p$.

$\mathbf{(SS)}$  For each $u\in S_p^{ss}$, one has that $u$ is unramified in $L/L'$, $L'_u=\Qp$ and $a_u = 1 + p - |\tilde{E}_u(\mathbb{F}_p)| = 0$, where $\tilde{E}_u$ is the reduction of $E$ at $u$.

\medskip
For any subset $R$ of $S$ and any extension $\mathcal{F}$ of $L'$, we shall write $R(\mathcal{F})$ for the set of primes of $\mathcal{F}$ above $R$.

Throughout, $F$ will always denote either $L$ or $L(\mu_p)$, where $\mu_p$ is the group of $p$th roots of unity. We then write $F^{\cyc}$ for the cyclotomic $\Zp$-extension of $F$ and $F_n$ the intermediate subfield of $F^{\cyc}$ with $|F_n:F|=p^n$.

\subsection{Some examples of elliptic curves with mixed reduction at $p$}

As mentioned in the introduction, the methods of Greenberg  \cite[Proposition 5.4]{G99} and Mazur \cite[Lemma 8.19]{Maz} can yield many examples of elliptic curves with mixed reduction type.
In this subsection, we give some examples of elliptic curves of our own with mixed good reduction at primes above $p$ which can be constructed easily. Consider the field $\Q(\sqrt{-1})$. Let $p$ be a prime $\equiv 5$ (mod $12$). Then there exist integers $a$ and $b$ such that $p = a^2+ b^2 = (a+b\sqrt{-1})(a-b\sqrt{-1})$, where $a+b\sqrt{-1}$ and $a-b\sqrt{-1}$ are primes in $\Z[\sqrt{-1}]$. Note that the integers $a$ and $b$ are not divisible by $p$.

Consider the elliptic curve $E: y^2 = x^3 + (a+b\sqrt{-1})x + a-b\sqrt{-1}$ which is defined over $\Q(\sqrt{-1})$.

\bl \label{examples}

The above elliptic curve is a global minimal Weierstrass equation with good supersingular reduction at the prime $a+b\sqrt{-1}$ and good ordinary reduction at the prime $a-b\sqrt{-1}$.
\el

\bpf
The $c_4$-invariant and discriminant of $E$ (in the sense of \cite[Chap. III, pp 42]{Sil})\footnote{Incidently, there is a misprint here. The $b_2$-invariant should be given by $a_1^2+4a_2$.} are given by
$-48(a+b\sqrt{-1})$ and
\[ -16\Big(4(a+b\sqrt{-1})^3 +27(a-b\sqrt{-1})^2\Big) \]
respectively. Since $c_4 = -2^4 \cdot3(a+b\sqrt{-1})$ and 3 is prime in $\Q(\sqrt{-1})$, it follows from \cite[Chap. VII, Remark 1.1]{Sil} that the defining equation of $E$ is minimal at all primes of $\Q(\sqrt{-1})$ except possibly $1+\sqrt{-1}$. Now note that $16 = (1+\sqrt{-1})^8$ is 12th-powerfree, and that the (Gaussian) norm of $4(a+b\sqrt{-1})^3 +27(a-b\sqrt{-1})^2$ is given by
$729p^2 + 16p^3 + 108(a^5 -10a^3b^2+5ab^4)$
which is an odd integer. This in turn implies that  $4(a+b\sqrt{-1})^3 +27(a-b\sqrt{-1})^2$ is not divisible by $1+\sqrt{-1}$. By \cite[Chap. VII, Remark 1.1]{Sil}, it then follows from this that the defining equation of $E$ is also minimal at $1+\sqrt{-1}$. In conclusion, the defining equation of $E$ is a global minimal Weierstrass equation.

Now, observe that $y^2 \equiv x^3 + 2a$ (mod $a+b\sqrt{-1}$). Since $p\equiv 5$ (mod $12$), we have
$p\equiv 1$ (mod $4$), the reduced curve is defined over $\mathbb{F}_p$. Since we also have $p\equiv 2$ (mod $3$), $(x^3+2a)^{(p-1)/2}$ has no terms of the form $x^{p-1}$. It then follows from \cite[Chap. V, Theorem 4.1(a)]{Sil} that the elliptic curve has good supersingular reduction at $a+b\sqrt{-1}$.

On the other hand,  we have $y^2 \equiv x^3 + 2ax$ (mod $a-b\sqrt{-1}$). Since $p$ splits completely in $\Q(\sqrt{-1})$, the reduced curve is
defined over $\mathbb{F}_p$  and so the coefficient of $x^{p-1}$ in $(x^3+2ax)^{(p-1)/2}$ is given by
\[{(p-1)/2\choose (p-1)/4}2^{(p-1)/4}a^{(p-1)/4}\]
which is nonzero in $\mathbb{F}_p$. Hence by \cite[Chap. V, Theorem 4.1(a)]{Sil}, the elliptic curve has good ordinary reduction at $a-b\sqrt{-1}$.
\epf

For instance, if $p=5$, the elliptic curve is given by $y^2 = x^3 + (2+\sqrt{-1})x + 2-\sqrt{-1}$ which has good supersingular reduction at $2+\sqrt{-1}$ and good ordinary reduction at $2-\sqrt{-1}$. By brute force, one can compute and factorize the discriminant of $E/\Q(\sqrt{-1})$ as
\[-16 ( 89 - 64 \sqrt{-1}) = -(1+\sqrt{-1})^8(6-5\sqrt{-1})(14+\sqrt{-1}), \]
where $6-5\sqrt{-1}$ is a prime of $\Q(\sqrt{-1})$ above $61$ and $14+\sqrt{-1}$ is a prime of $\Q(\sqrt{-1})$ above $197$. We also mention that this curve has Mordell rank 2 (over $\Q(\sqrt{-1})$) by Sage.

If $p=17$, we have $y^2 = x^3 + (4+\sqrt{-1})x + 4-\sqrt{-1}$ which has good supersingular reduction at $4+\sqrt{-1}$ and good ordinary reduction at $4-\sqrt{-1}$. One can compute and factorize the discriminant of $E/\Q(\sqrt{-1})$ to be
\[-16 (613  -28 \sqrt{-1}) = -(1+\sqrt{-1})^8(6+5\sqrt{-1})(58-53\sqrt{-1}), \]
where $6+5\sqrt{-1}$ is a prime of $\Q(\sqrt{-1})$ above $61$ and $58-53\sqrt{-1}$ is a prime of $\Q(\sqrt{-1})$ above $6173$. We note that this curve has Mordell rank 1 (over $\Q(\sqrt{-1})$) by Sage.

Unfortunately, we do not know how to compute the Iwasawa invariants of the mixed signed Selmer group for any of these examples at this point of writing.

\subsection{Signed Selmer groups} \label{subsec signed Selmer}

In this subsection, we give the definition of the signed Selmer groups. The notion of the signed Selmer groups was first conceived by Kobayashi \cite{Kob} to handle elliptic curves with good supersingular reduction at all primes above $p$. Later, Kim \cite{Kim14} introduced the signed Selmer groups which allowed mixed signs. Recently, Kitajima and Otsuki \cite{KO} considered the signed Selmer groups attached to elliptic curves with mixed good reduction types at primes above $p$.
Mimicking these prior works, the authors of this paper introduced the signed Selmer groups which allow mixed signs and mixed good reduction types (see \cite{AL}). In this paper, we will slightly refine this even further by allowing the presence of multiplicative reduction primes above $p$ (also see \cite{LLAk}).

We begin by taking care of the local conditions at the supersingular primes.
By $\mathbf{(SS)}$, every prime in $S_p^{ss}(F)$ is totally ramified in $F^{\cyc}/F$. In particular, for each such prime $v$, there is a unique prime of $F_n$ lying above the said prime which, by abuse of notation, is still denoted by $v$. We then write $\hat{E}_v$ for the formal group $E$ over $L'_u$, where $u$ is a prime of $L'$ below $v$.

For each $\overrightarrow{s}=(s_v)_{v\in S_p^{ss}(F)}\in\{\pm\}^{S_p^{ss}(F)}$, we write \[\mathcal{H}^{\overrightarrow{s}}_n = \bigoplus_{v\in S_p^{ss}(F)}\frac{H^1(F_{n,v},E(p))}{\hat{E}_v^{s_v}(F_{n,v})\ot\Qp/\Zp}.\]
The signed Selmer groups are then defined by
\[\Sel^{\overrightarrow{s}}(E/F_n) = \ker \left(H^1(G_S(F_n),E(p))\stackrel{\psi^{\overrightarrow{s}}}{\lra} \mathcal{H}^{\overrightarrow{s}}_n\times\bigoplus_{w\in S_p^{o}(F_n)}\frac{H^1(F_{n,w},E(p))}{E(F_{n,w})\ot\Qp/\Zp}\times\bigoplus_{w\in S'(F_n)}H^1(F_{n,w},E(p)) \right).\]

Note that in the event that $S_p^{ss}=\emptyset$, the above definition coincides with the usual $p$-primary Selmer group $\Sel(E/F_n)$. In the event that all the signs occuring in the definition are $-$ signs, the signed Selmer group is then denoted by $\Sel^{\overrightarrow{-}}(E/F_n)$. We also note that one always has $\Sel^{\overrightarrow{s}}(E/F)= \Sel(E/F)$ regardless of the presence of supersingular primes. For a general $n$, the signed Selmer group and the classical Selmer group fit into the following commutative diagram
\[   \xymatrixrowsep{0.25in}
\xymatrixcolsep{0.15in}\entrymodifiers={!! <0pt, .8ex>+} \SelectTips{eu}{}\xymatrix{
    0 \ar[r]^{} & \Sel^{\overrightarrow{s}}(E/F_n) \ar[d]_{\al} \ar[r] &  H^1(G_S(F_n), E(p))
    \ar@{=}[d] \ar[r]^(.27){\psi^{\overrightarrow{s}}} & \displaystyle\mathcal{H}^{\overrightarrow{s}}_n\times\bigoplus_{w\in S_p^{o}(F_n)}\frac{H^1(F_{n,w},E(p))}{E(F_{n,w})\ot\Qp/\Zp}\times\bigoplus_{w\in S'(F_n)}H^1(F_{n,w},E(p))\ar[d]\\
    0 \ar[r]^{} & \Sel(E/F_{n}) \ar[r]^{} & H^1(G_S(F_n), E(p)) \ar[r]^(.3){\phi} & \displaystyle\bigoplus_{w\in S_p(F_n)}\frac{H^1(F_{n,w},E(p))}{E(F_{n,w})\ot\Qp/\Zp}\times\bigoplus_{w\in S'(F_n)}H^1(F_{n,w},E(p)) } \]
with exact rows. Denote by $\psi^{\overrightarrow{s}}_{ss}$ the map from $\Sel(E/F_n)$ to $\mathcal{H}^{\overrightarrow{s}}_n$ that is induced by $\psi^{\overrightarrow{s}}$. Then one has the following equivalent description of the signed Selmer groups.

\bl \label{PM Selmer} We have the following identification
\[\Sel^{\overrightarrow{s}}(E/F_n) = \ker \left(\Sel(E/F_n)\stackrel{\psi^{\overrightarrow{s}}_{ss}}{\lra} \mathcal{H}^{\overrightarrow{s}}_n\right).\]
\el

\bpf
This follows from a straightforward analysis of the above commutative diagram with the definition of the signed Selmer groups.
\epf

We then define $\Sel^{\overrightarrow{s}}(E/F^{\cyc})=\ilim_n \Sel^{\overrightarrow{s}}(E/F_n)$ and $\mathcal{H}^{\overrightarrow{s}}_{\infty}= \ilim_n \mathcal{H}^{\overrightarrow{s}}_n$.
It is not difficult to verify that  $\Sel^{\overrightarrow{s}}(E/F^{\cyc})$ is cofinitely generated over $\Zp\ps{\Ga}$, where $\Ga=\Gal(F^{\cyc}/F)$. In fact, one expects the following conjecture which is a natural extension of Mazur \cite{Maz} and Kobayashi \cite{Kob}.

\medskip \noindent \textbf{Conjecture.} $\Sel^{\overrightarrow{s}}(E/F^{\cyc})^{\vee}$ is a torsion $\Zp\ps{\Ga}$-module.

\medskip
When $E$ has good ordinary reduction at all primes above $p$, the above conjecture is precisely Mazur's conjecture \cite{Maz} which is known in the case when $E$ is defined over $\Q$ and $F$ an abelian extension of $\Q$ (see \cite{K}). For an elliptic curve over $\Q$ with good supersingular reduction at $p$, this conjecture was established by Kobayashi
(cf. \cite{Kob}; also see \cite{BL} for some recent progress on this conjecture).

For our subsequent discussion, we require another description of $\Sel^{\overrightarrow{s}}(E/F^{\cyc})$.
Let $v$ be a prime of $F$ above $p$ which is not a good supersingular reduction prime of $E$. Then from Subsection \ref{ordinary subsection}, we have the following short exact sequence
\[ 0\lra C_v \lra E(p)\lra D_v\lra 0\]
of $\Gal(\bar{F}_v/F_v)$-modules, where both $C_v$ and $D_v$ are cofree with $\Zp$-corank 1. If $w$ is a prime of $F^{\cyc}$ above $v$, we sometimes write $D_w$ for $D_v$. Now the results of Coates and Greenberg \cite[Propositions 4.3 and 4.8]{CG} tell us that $H^1(F^{\cyc}_{w},E)(p) \cong H^1(F^{\cyc}_{w}, D_v)$ which is the key to obtaining the following alternative useful description of the signed Selmer groups.

\bp \label{selmer alternative} We have the following two equivalent descriptions of the signed Selmer group, namely \[\Sel^{\overrightarrow{s}}(E/F^{\cyc})\cong \ker\left( H^1(G_S(F^{\cyc}),E(p))\stackrel{\psi^{\overrightarrow{s}}}{\lra} \mathcal{H}^{\overrightarrow{s}}_{\infty}\times\bigoplus_{w\in S_p^{o}(F^{\cyc})} H^1(F^{\cyc}_w, D_w)\times\bigoplus_{w\in S'(F^{\cyc})}H^1(F^{\cyc}_{w},E(p))\right)\]
\[ \cong \ker\left( H^1(G_S(F^{\cyc}),E(p))\stackrel{\psi^{\overrightarrow{s}}}{\lra} \mathcal{H}^{\overrightarrow{s}}_{\infty}\times\bigoplus_{w\in S_p^{o}(F^{\cyc})} H^1_w(E/F^{\cyc})\times\bigoplus_{w\in S'(F^{\cyc})}H^1(F^{\cyc}_{w},E(p))\right),\]
where
\[H^1_w(E/F^{\cyc}) = \begin{cases}  H^1(F^{\cyc, ur}_{w}, D_w),& \mbox{if $w$ is a prime of good ordinary reduction} \\
& \mbox{or non-split multiplicative reduction}, \\
      H^1(F^{\cyc}_w, \Qp/\Zp), & \mbox{if $w$ is a prime of split multiplicative reduction.} \end{cases}
 \]
\ep

\bpf
The first isomorphism follows immediately from the discussion before the proposition. If $D_v(F^{\cyc}_w)$ is infinite, then it is $\Qp/\Zp$ as $\Gal(\overline{F^{\cyc}_w}/F^{\cyc}_w)$-module by Lemma \ref{ordinary}. The same lemma says that if $w$ is a prime of good ordinary reduction or non-split multiplicative reduction, then $H^1(F^{\cyc}_{w}, D_v)\cong H^1(F^{\cyc, ur}_{w}, D_v)$. Hence we may replace $H^1(F^{\cyc}_{w}, D_w)$ by $H^1_w(E/F^{\cyc})$ for $w\in S_p^o(F^{\cyc})$ in the first isomorphism and still recover $\Sel^{\overrightarrow{s}}(E/F^{\cyc})$. (Note: in view of this identification, the localisation maps in both isomorphisms are then identified with the map $\psi^{\overrightarrow{s}}$ as given in the definition of the signed Selmer groups, and by abuse of notation, we also denote both localisation maps by $\psi^{\overrightarrow{s}}$ as stated in the proposition.)
\epf

We shall require another equivalent description of the signed Selmer groups on the level of $F$ and $F^{\cyc}$. For each $n$, we define the strict signed Selmer group
\[\Sel^{\overrightarrow{s},str}(E/F_n)= \ker\left( H^1(G_S(F_n),E(p))\lra \bigoplus_{w\in S(F_n)}\displaystyle\frac{H^1(F_{n,w},E(p))}{L_w(E/F_n)} \right),\]
where
\[ L_w(E/F_n) = \begin{cases}  \Big(E^{s_w}(F^{\cyc}_{w})\ot\Qp/\Zp\Big)^{\Ga_n},& \mbox{if $w\in S_p^{ss}(F_n)$}, \\
 \Big(\im\left( H^1(F_{n,w},C_w)\lra H^1(F_{n,w}, E(p))\right)\Big)_{div},& \mbox{if $w\in S_p^{o}(F_n)$}, \\
 0, & \mbox{if $w\in S'(F_n)$,} \end{cases}\]
where $\Ga_n$ denotes $\Gal(F^{\cyc}/F_n)$ and $(M)_{div}$ is the maximal $p$-divisible subgroup of $M$.

\br
Note that in the case when $S_p^{ss}=\emptyset$, this is the strict Selmer group as defined in the sense of Greenberg \cite{G89} (also see \cite{Guo}). Our choice of naming the above as the strict signed Selmer group is inspired by this observation.
\er

We now give a result which compares the strict signed Selmer groups with the usual $p$-primary Selmer group at the level of $F$, and compares the strict signed Selmer groups with the signed Selmer groups at the level of $F^{\cyc}$. We should mention that the strict Selmer groups $\Sel^{\overrightarrow{s},str}(E/F_n)$ need not coincide with $\Sel^{\overrightarrow{s}}(E/F_n)$ in general. Before stating and proving our result, we introduce another assumption that will play a role in subsequent discussion especially in the presence of a $``+"$ sign in the signed Selmer groups.

\medskip
$\mathbf{(S+)}$ For every $v\in S_p^{ss}(L)$, assume further that $[L_v:\Qp]$ is not divisible by 4.

\bp \label{strict selmer=selmer}
 The following statements are valid.
\begin{enumerate}
\item[$(a)$] $\Sel^{\overrightarrow{s}}(E/F^{\cyc})\cong \ilim_n \Sel^{\overrightarrow{s},str}(E/F_n)$ for every $\overrightarrow{s}$.
\item[$(b)$] $\Sel^{\overrightarrow{-},str}(E/F) = \Sel(E/F)$.
\item[$(c)$] If one assumes further that $\mathbf{(S+)}$ holds, then we also have $\Sel^{\overrightarrow{s},str}(E/F) = \Sel(E/F)$ for every $\overrightarrow{s}$.
\end{enumerate}\ep

\bpf
To prove assertions (b) and (c), it suffices to show that $L_v(E/F) = E(F_v)\ot\Qp/\Zp$ for all $v\in S$. This is clear for $v\in S'(F)$. For $v\in S_p^o(F)$, this is a result of Coates-Greenberg (cf. \cite[Proposition 4.5]{CG}). Now suppose that $v\in S^{ss}_p(F)$. We have a natural inclusion
\[ E(F_v)\ot \Qp/\Zp \hookrightarrow \big(E^{\pm}(F^{\cyc}_w)\ot \Qp/\Zp\big)^{\Ga}.\]
By $\mathbf{(S+)}$ and Proposition \ref{Kita-O results} (or \cite[Corollary 3.25]{KO}), we have that $\big(E^{\pm}(F^{\cyc}_v)\ot \Qp/\Zp\big)^{\Ga}$ is a cofree $\Zp$-module of $\Zp$-corank $[F_v:\Qp]$. On the other hand, it follows from Mattuck's theorem \cite{Mat} that $E(F_v)\ot \Qp/\Zp$ is also a cofree $\Zp$-module of $\Zp$-corank $[F_v:\Qp]$. Therefore, the inclusion has to be an isomorphism.

For assertion (a), we note that $L_w(E/F_n) = E(F_{n,w})\ot\Qp/\Zp$ for all $w\in S_p^o(F_n)$ and so the local condition at $S_p^o(F_n)$ for the strict signed Selmer group and the signed Selmer group at each $F_n$ (and hence $F^{\cyc}$)) are the same. Similarly, one has the same conclusion for primes in $S'_p(F^{\cyc})$. Let $v$ be a prime $S_p^{ss}$. Recall that we also write $v$ for the prime of $F^{\cyc}$ above $v$. In general,
$\Big(E^{s_v}(F^{\cyc}_{v})\ot\Qp/\Zp\Big)^{\Ga_n}$ and $E^{s_v}(F_{n,v})\ot\Qp/\Zp$ need not agree (one can see this by comparing their $\Zp$-coranks). But upon taking limit, they are both equal to
$E^{s_v}(F^{\cyc}_{v})\ot\Qp/\Zp$. In conclusion, we have established assertion (a).
\epf

\subsection{Torsionness of signed Selmer groups}

As noted in the previous subsection, $\Sel^{\overrightarrow{s}}(E/F^{\cyc})$ is expected to be cotorsion. In this section, we consider an equivalent characterization of this property which is precisely the content of the next proposition.

\bp \label{torsion surjective H2}
$\Sel^{\overrightarrow{s}}(E/F^{\cyc})$ is a cotorsion $\Zp\ps{\Ga}$-module if and only if $H^2(G_S(F^{\cyc}),E(p))=0$ and $\psi^{\overrightarrow{s}}$ is surjective.
\ep

The above result is well-known for the usual $p$-primary Selmer group when the elliptic curve $E$ has either good ordinary reduction or multiplicative reduction at all primes above $p$ (for instances, see \cite[Proposition 7.2]{HV} or \cite[Proposition 3.3]{LimMHG}). When the elliptic curve $E$ has good supersingular reduction at all primes above $p$, this was
somewhat partially established in \cite[Proposition 2.4]{Kim09} (also see \cite[Proposition 3.10]{KimPM}).
That the analogue statement also holds in this mixed reduction situation was made aware to us by Antonio Lei and Ramdorai Sujatha, and we thank them for this and for sharing their work \cite{LeiSu} on the subject.

As noted above, if $S_p^{ss}=\emptyset$, Proposition \ref{torsion surjective H2} is then a consequence of \cite[Proposition 3.3]{LimMHG}. Hence it remains to establish the equivalence as asserted in Proposition \ref{torsion surjective H2} under the assumption that $S_p^{ss}\neq\emptyset$. As a start, we record a few preparatory lemmas.

\bl \label{E(p)=0} Suppose that $S_p^{ss}\neq\emptyset$. Then $E(F^{\cyc})(p)=0$.
\el

\bpf
 For $w\in S_p^{ss}(F)$, it follows from Lemma \ref{supersingular points} that $E(E^{\cyc}_w)(p)=0$. Since $E(F^{\cyc})(p)$ is contained in $E(E^{\cyc}_w)(p)$, the conclusion of the lemma follows.
\epf

For the subsequent discussion, we shall write $\Hi^i(F^{\cyc}/F, T_pE) = \plim_n H^i(G_S(F_n), T_pE)$, where $T_pE$ is the Tate module of the elliptic curve $E$.

\bl \label{H1Iw} Suppose that $S_p^{ss}\neq\emptyset$. Then $\Hi^1(F^{\cyc}/F, T_pE)$ is a torsion free $\Zp\ps{\Ga}$-module.
\el

\bpf
By considering the low degree terms of the spectral sequence of Jannsen
\[ \Ext^i_{\Zp\ps{\Ga}}\big(H^j(G_S(F^{\cyc}),E(p))^{\vee},\Zp\ps{\Ga}\big) \Longrightarrow \Hi^{i+j}(F^{\cyc}/F, T_pE)\]
(cf. \cite[Theorem 1]{Jan}), we obtain the following exact sequence
\[ 0\lra \Ext^1_{\Zp\ps{\Ga}}\big(E(F^{\cyc})(p))^{\vee},\Zp\ps{\Ga}\big) \lra \Hi^1(F^{\cyc}/F, T_pE) \lra \Ext^0_{\Zp\ps{\Ga}}\big(H^1(G_S(F^{\cyc}),E(p))^{\vee},\Zp\ps{\Ga}\big). \]
In view of Lemma \ref{E(p)=0}, the leftmost term vanishes. It then follows that $\Hi^1(F^{\cyc}/F, T_pE)$ injects into an $\Ext^0$-term which is a reflexive $\Zp\ps{\Ga}$-module by \cite[Corollary 5.1.3]{NSW}. Therefore, $\Hi^1(F^{\cyc}/F, T_pE)$ must be torsion free. \epf

We can now give a proof of Proposition \ref{torsion surjective H2} (compare with \cite[Proposition 4.4]{LeiSu} and \cite[Proposition 3.3]{LimMHG}).

\bpf[Proof of Proposition \ref{torsion surjective H2}]
As noted above, it suffices to prove the proposition under the assumption that $S_p^{ss}\neq\emptyset$ which we do. For the proof, we shall make use of the first description of $ \Sel^{\overrightarrow{s}}(E/F^{\cyc})$ in Proposition \ref{selmer alternative}.
By \cite[Proposition A.3.2]{PR00}, there is an exact sequence
\[0\lra \Sel^{\overrightarrow{s}}(E/F^{\cyc})\lra H^1(G_S(F^{\cyc}),E(p))\stackrel{\psi^{\overrightarrow{s}}}{\lra} \mathcal{H}^{\overrightarrow{s}}_{\infty}\times\bigoplus_{w\in S_p^{o}(F^{\cyc})}H^1(F^{\cyc}_{w},D_w)\times\bigoplus_{w\in S'(F^{\cyc})}H^1(F^{\cyc}_{w},E(p)) \]
\[ \lra \mathfrak{S}^{\overrightarrow{s}}(E/F^{\cyc})^{\vee}\lra H^2(G_S(F^{\cyc}),E(p))\lra 0,\]
where $\mathfrak{S}^{\overrightarrow{s}}(E/F^{\cyc})$ is a $\Zp\ps{\Ga}$-submodule of $\Hi^1(F^{\cyc}/F, T_pE)$. (For the precise definition of $\mathfrak{S}^{\overrightarrow{s}}(E/F^{\cyc})$, we refer readers to loc. cit. For our purposes, the submodule theoretical information suffices.) Standard corank calculations \cite[Propositions 1-3]{G89} tell us that
\[ \corank_{\Zp\ps{\Ga}}\big(H^1(G_S(F^{\cyc}),E(p))\big)- \corank_{\Zp\ps{\Ga}}\big(H^2(G_S(F^{\cyc}),E(p))\big) =[F: \Q], \]
 \[\corank_{\Zp\ps{\Ga}}\left(\bigoplus_{w\in S_p^{o}(F^{\cyc})}H^1(F^{\cyc}_{w},D_w)\right) = \sum_{v\in S_p^{o}(F)}[F_v:\Qp]\] and
\[\corank_{\Zp\ps{\Ga}}\left(\bigoplus_{w\in S'(F^{\cyc})}H^1(F^{\cyc}_{w},E(p))\right) = 0.\]
On the other hand, it follows from \cite[Proposition 3.32]{KO} that \[\corank_{\Zp\ps{\Ga}}\big(\mathcal{H}^{\overrightarrow{s}}_{\infty}\big) = \sum_{v\in S_p^{ss}(F)}[F_v:\Qp].\]
It is now clear from these formulas and the above exact sequence that $\Sel^{\overrightarrow{s}}(E/F^{\cyc})$ is a cotorsion $\Zp\ps{\Ga}$-module if and only if $\mathfrak{S}^{\overrightarrow{s}}(E/F^{\cyc})$ is a torsion $\Zp\ps{\Ga}$-module. Since $\mathfrak{S}^{\overrightarrow{s}}(E/F^{\cyc})$ is contained in $\Hi^1(F^{\cyc}/F, T_pE)$ which is torsion free by Lemma \ref{H1Iw}, the latter statement holds if and only if $\mathfrak{S}^{\overrightarrow{s}}(E/F^{\cyc})=0$. In view of the exact sequence in the beginning of the proof, this is precisely equivalent to saying that
$H^2(G_S(F^{\cyc}),E(p))=0$ and $\psi^{\overrightarrow{s}}$ is surjective. Hence we have established the proposition.
\epf

\subsection{Non-primitive signed Selmer group}

In comparing Selmer groups of two congruent elliptic curves, it is a standard procedure to do so via an appropriate non-primitive variant of the Selmer group (see \cite{GV, Kim09}). We shall follow this strategy and begin by introducing the non-primitive variant of the signed Selmer group. This in turn is denoted and defined by
\[\Sel^{\overrightarrow{s}}_{non}(E/F^{\cyc}) = \ker \left(H^1(G_S(F^{\cyc}),E(p))\stackrel{\psi^{\overrightarrow{s}}_{non}}\lra \mathcal{H}^{\overrightarrow{s}}_{\infty}\times\bigoplus_{w\in S_p^{o}(F^{\cyc})}H^1_w(E/F^{\cyc}) \right).\]

 The non-primitive signed Selmer group and the signed Selmer group are connected via the following commutative diagram
\[   \xymatrixrowsep{0.25in}
\xymatrixcolsep{0.15in}\entrymodifiers={!! <0pt, .8ex>+} \SelectTips{eu}{}\xymatrix{
    0 \ar[r]^{} & \Sel^{\overrightarrow{s}}(E/F^{\cyc}) \ar[d] \ar[r] &  H^1(G_S(F^{\cyc}), E(p))
    \ar@{=}[d] \ar[r]^(.3){\psi^{\overrightarrow{s}}} & \displaystyle\mathcal{H}^{\overrightarrow{s}}_{\infty}\times\bigoplus_{w\in S_p^{o}(F^{\cyc})}H^1_w(E/F^{\cyc})\times\bigoplus_{w\in S'(F^{\cyc})}H^1(F^{\cyc}_{w},E(p))\ar[d]\\
    0 \ar[r]^{} & \Sel^{\overrightarrow{s}}_{non}(E/F^{\cyc}) \ar[r]^{} & H^1(G_S(F^{\cyc}), E(p)) \ar[r]^(.4){\psi^{\overrightarrow{s}}_{non}} &\displaystyle\mathcal{H}^{\overrightarrow{s}}_{\infty}\times\bigoplus_{w\in S_p^{o}(F^{\cyc})}H^1_w(E/F^{\cyc}) } \]
with exact rows. An application of the snake lemma then gives following exact sequence
\[0\lra \Sel^{\overrightarrow{s}}(E/F^{\cyc}) \lra \Sel^{\overrightarrow{s}}_{non}(E/F^{\cyc})\lra \bigoplus_{w\in S'(F^{\cyc})}H^1(F^{\cyc}_{w},E(p)). \]
We can now record an important result on the structure of the non-primitive signed Selmer group.

\bp \label{no finite submodule}
Suppose that $\Sel^{\overrightarrow{s}}(E/F^{\cyc})$ is cofinitely generated over $\Zp$. Then the following statements hold.

$(a)$ $\Sel^{\overrightarrow{s}}_{non}(E/F^{\cyc})^{\vee}$ has no nonzero finite $\Zp\ps{\Ga}$-submodules.

$(b)$ There is a short exact sequence
\[0\lra \Sel^{\overrightarrow{s}}(E/F^{\cyc}) \lra \Sel^{\overrightarrow{s}}_{non}(E/F^{\cyc})\lra \bigoplus_{w\in S'(F^{\cyc})}H^1(F^{\cyc}_{w},E(p)) \lra 0. \]
\ep

\bpf
 Since $\Sel^{\overrightarrow{s}}(E/F^{\cyc})$ is cofinitely generated over $\Zp$, it is also a cotorsion $\Zp\ps{\Ga}$-module. It then follows from Proposition \ref{torsion surjective H2} that $H^2(G_S(F^{\cyc}),E(p))=0$ and $\psi^{\overrightarrow{s}}$ is surjective. The former implies that $H^1(G_S(F^{\cyc}),E(p))^{\vee}$ has no nonzero finite $\Zp\ps{\Ga}$-submodules (cf. \cite[Proposition 5]{G89}), whereas the latter, when combined with a diagram-chasing argument, yields the short exact sequence in (b) and the surjectivity of $\psi^{\overrightarrow{s}}_{non}$.

 Now if $S_p^{ss}$ is non-empty, it follows from Lemma \ref{E(p)=0} that $E(F^\cyc)[p]=0$. We can then follow the approach of \cite[Proposition 4.14]{G99} and \cite[Theorem 3.14]{KimPM} to show that $\Sel^{\overrightarrow{s}}_{non}(E/F^{\cyc})^{\vee}$ has no nonzero finite $\Zp\ps{\Ga}$-submodules. In the event that $S_p^{ss}$ is empty, the approach of Greenberg \cite[Proposition 4.14]{G99} cannot carry over, since we are not assuming that $E(F)[p]=0$. Fortunately, in this situation, we can use the approach in \cite[Proposition 2.8]{LimSu}.
\epf

We record the following corollary of the preceding proposition.

\bc \label{Zp-rank} The following statements are equivalent.
\begin{enumerate}
\item[$(a)$] $\Sel^{\overrightarrow{s}}(E/F^{\cyc})$ is a cofinitely generated $\Zp$-module.
\item[$(b)$] $\Sel^{\overrightarrow{s}}_{non}(E/F^{\cyc})$ is a cofinitely generated $\Zp$-module.
\item[$(c)$] $\Sel^{\overrightarrow{s}}_{non}(E/F^{\cyc})[p]$ is finite.
     \end{enumerate}
Moreover, in the event of that one of (and hence all of) the above statements holds, one has
\[\ba{rl}\dim_{\mathbb{F}_p}\Big(\Sel^{\overrightarrow{s}}_{non}(E/F^{\cyc})[p]\Big)\hspace{-0.12in}&= \corank_{\Zp}\Big(\Sel^{\overrightarrow{s}}_{non}(E/F^{\cyc})\Big) \\
&=  \corank_{\Zp}\Big(\Sel^{\overrightarrow{s}}(E/F^{\cyc})\Big) + \displaystyle\sum_{w\in S'(F^{\cyc})}H^1(F^{\cyc}_{w},E(p)).  \ea \]
\ec

\bpf
 The equivalence of (a) and (b) is an immediate consequence of Proposition \ref{no finite submodule}(b). The equivalence of (b) and (c) is a consequence of Nakayama lemma. Now if $\Sel^{\overrightarrow{s}}_{non}(E/F^{\cyc})$ is cofinitely generated over $\Zp$, then
 $\Sel^{\overrightarrow{s}}_{non}(E/F^{\cyc})/p$ is a finite quotient of $\Sel^{\overrightarrow{s}}_{non}(E/F^{\cyc})$, and by Proposition \ref{no finite submodule}(a), it has to be trivial. The first equality now follows from this. The final equality follows from Proposition \ref{no finite submodule}(b) again.
 \epf

We now consider the non-primitive mod-$p$ signed Selmer group $\Sel^{\overrightarrow{s}}_{non}(E[p]/F^{\cyc})$ which is defined to be
\[ \ker \left(H^1(G_S(F^{\cyc}),E[p])\lra\mathcal{I}_{\infty}^{ss}\times\bigoplus_{w\in S_p^{0}(F^{\cyc})}H^1_w(E[p]/F^{\cyc}) \right),\]
where $\displaystyle\mathcal{I}_{\infty}^{ss} =\bigoplus_{w\in S_p^{ss}(F^{cyc})}\frac{H^1(F^{\cyc}_{w_i},E[p])}{E^{s_i}(F^{\cyc}_{w})/p}$ and
\[H^1_w(E[p]/F^{\cyc}) = \begin{cases}  H^1(F^{\cyc, ur}_{w},D_w[p]),& \mbox{if $D_v(F^{\cyc}_w)$ is finite}, \\
      H^1(F^{\cyc}_w, D_w[p]), & \mbox{if $D_v(F^{\cyc}_w)$ is infinite.} \end{cases}
 \]
This fits into the following commutative diagram
\[   \xymatrixrowsep{0.25in}
\xymatrixcolsep{0.15in}\entrymodifiers={!! <0pt, .8ex>+} \SelectTips{eu}{}\xymatrix{
    0 \ar[r]^{} & \Sel^{\overrightarrow{s}}_{non}(E[p]/F^{\cyc}) \ar[d] \ar[r] &  H^1(G_S(F^{\cyc}), E[p])
    \ar[d]^{h} \ar[r]^(){} & \displaystyle\mathcal{I}_{\infty}^{ss}\times\bigoplus_{w\in S_p^{o}(F^{\cyc})}H^1_w(E[p]/F^{\cyc})\ar[d]^(.6){\oplus c_w}\\
    0 \ar[r]^{} & \Sel^{\overrightarrow{s}}_{non}(E/F^{\cyc})[p] \ar[r]^{} & H^1(G_S(F^{\cyc}), E(p))[p] \ar[r] &\displaystyle\mathcal{H}_{\infty}^{ss}[p]\times\bigoplus_{w\in S_p^{o}(F^{\cyc})}H^1_w(E/F^{\cyc})[p] } \]
with exact rows.

\bl \label{mod p} $\Sel^{\overrightarrow{s}}_{non}(E/F^{\cyc})$ is cofinitely generated over $\Zp$ if and only if $\Sel^{\overrightarrow{s}}_{non}(E[p]/F^{\cyc})$ is finite. Furthermore, in the event of such, we have the following equality
\[\dim_{\mathbb{F}_p}\Big(\Sel^{\overrightarrow{s}}_{non}(E[p]/F^{\cyc})\Big)= \corank_{\Zp}\Big(\Sel^{\overrightarrow{s}}_{non}(E/F^{\cyc})\Big) + \dim_{\mathbb{F}_p}\Big((E(F^{\cyc})[p]\Big).  \]
\el

\bpf
The map $h$ in the above commutative diagram is surjective with $\ker h =E(F^{\cyc})(p)/p$. Since $E(F^{\cyc})(p)$ is finite by \cite{Ri}, we have $\dim_{\mathbb{F}_p}\big((E(F^{\cyc})[p]\big) = \dim_{\mathbb{F}_p}\big(\ker h\big)$. We now show that each $c_w$ is an injection. Suppose that $w$ lies above $v\in S_p^{o}(F)$. Then one has
\[\ker c_w = \begin{cases} D_w(F^{\cyc, ur})/p,& \mbox{if $D_v(F^{\cyc}_w)$ is finite}, \\
     D_w(F^{\cyc})/p, & \mbox{if $D_v(F^{\cyc}_w)$ is infinite.} \end{cases}
 \]
 But in either of the two cases,  $\ker c_w$ is a mod-$p$ quotient of a divisible group and hence must be zero.

Suppose that $w\in S_p^{ss}(F^{\cyc})$. Then we have the following short exact sequence
\[0\lra E^{\pm}(F^{\cyc}_w)\ot \Qp/\Zp \lra H^1(F^{\cyc}_w, E(p))\lra \frac{H^1(F^{\cyc}_w, E(p))}{E^{\pm}(F^{\cyc}_w)\ot \Qp/\Zp}\lra 0. \]
Since $E^{\pm}(F^{\cyc}_w)\ot \Qp/\Zp$ is $p$-divisible, the above exact sequence yields the following short exact sequence
\[0\lra E^{\pm}(F^{\cyc}_w)/p \lra H^1(F^{\cyc}_w, E(p))[p]\lra \left(\frac{H^1(F^{\cyc}_w, E(p))}{E^{\pm}(F^{\cyc}_w)\ot \Qp/\Zp}\right)[p]\lra 0. \]
Since $E(F^{\cyc}_w)(p)=0$, the term $H^1(F^{\cyc}_w, E(p))[p]$ identifies with $H^1(F^{\cyc}_w, E[p])$. It then follows from this that $c_w$ is also injective in this case.

In conclusion, the above argument shows that $\Sel^{\overrightarrow{s}}_{non}(E/F^{\cyc})[p]$ is finite if and only if $\Sel^{\overrightarrow{s}}_{non}(E[p]/F^{\cyc})$ is finite. Furthermore, in the event that this finiteness property holds, the argument also shows that
\[\ba{rl}\dim_{\mathbb{F}_p}\Big(\Sel^{\overrightarrow{s}}_{non}(E[p]/F^{\cyc})\Big)\hspace{-0.12in}&= \dim_{\mathbb{F}_p}\Big(\Sel^{\overrightarrow{s}}_{non}(E/F^{\cyc})[p]\Big) + \dim_{\mathbb{F}_p}\Big((E(F^{\cyc})[p]\Big).  \ea \]
The conclusions of the lemma will then follow from combining these observations with Corollary \ref{Zp-rank}.
\epf

\br
Note that if $S_p^{ss}\neq \emptyset$, then we even have $\ker h=0$ and that \[\Sel^{\overrightarrow{s}}_{non}(E[p]/F^{\cyc})\cong \Sel^{\overrightarrow{s}}_{non}(E/F^{\cyc})[p].\]
\er

\bp \label{mod p rank}
$\Sel^{\overrightarrow{s}}(E/F^{\cyc})$ is cofinitely generated over $\Zp$ if and only if $\Sel^{\overrightarrow{s}}_{non}(E[p]/F^{\cyc})$ is finite. Moreover, if this is so, we have
\[ \corank_{\Zp}\Big(\Sel^{\overrightarrow{s}}(E/F^{\cyc})\Big) + \sum_{w\in S'(F^{\cyc})}\corank_{\Zp}\Big(H^1(F^{\cyc}_w,E(p))\Big) + \dim_{\mathbb{F}_p}\Big((E(F^{\cyc})[p]\Big)\] \[= \dim_{\mathbb{F}_p}\Big(\Sel^{\overrightarrow{s}}_{non}(E[p]/F^{\cyc})\Big). \]
\ep

\bpf
The first assertion is clear from the exact sequence before Proposition \ref{no finite submodule}, whereas the second is a consequence of Corollary \ref{Zp-rank} and Lemma \ref{mod p}.
\epf

\section{Main results} \label{main results}

In this section, we shall prove our main results. Let $E_1$ and $E_2$ be two elliptic curves defined over $L'$ with $E_1[p]\cong E_2[p]$ as $\Gal(\overline{L'}/L)$-modules, where $L$ is a finite extension of $L'$. We shall retain the notation and settings from the previous section. In particular, the two elliptic curves are assumed to satisfy $\mathbf{(NA)}$ and $\mathbf{(SS)}$. As mentioned in the introduction, we need to be able to transfer this global information into local cohomology groups at primes above $p$ and this leads us to impose the following two hypotheses.

\medskip
$\mathbf{(RED)}$ The elliptic curves $E_1$ and $E_2$ have the same reduction type for each prime of $L'$ above $p$. Furthermore, in the event that $E_1$ and $E_2$ have non-split multiplicative reduction at a prime $u$ of $L'$ above $p$, suppose further $E_1$ and $E_2$ have the same reduction type at all primes of $F$ above $u$.

$\mathbf{(RU)}$ For each $u\in S_p^{o}$ such that $u$ is a prime of good ordinary reduction or non-split multiplicative reduction, suppose that $L'_u$ does not contain $\mu_p$.

\subsection{Comparing ranks of Selmer groups}

We can now state and prove the following result which generalizes that of Greenberg-Vatsal \cite{GV} and Kim \cite{Kim09}.

\bt \label{congruent theorem}
Let $E_1$ and $E_2$ be two elliptic curves defined over $L'$ with $E_1[p]\cong E_2[p]$ as $\Gal(\overline{L'}/L)$-modules. Let $L$ be a finite extension of $L'$. Write $F$ for either $L$ or $L(\mu_p)$. Suppose that $\mathbf{(NA)}$, $\mathbf{(SS)}$, $\mathbf{(RED)}$ and $\mathbf{(RU)}$ hold for $E_1$ and $E_2$.

Then $\Sel^{\overrightarrow{s}}(E_1/F^{\cyc})$ is cofinitely generated over $\Zp$ if and only if
$\Sel^{\overrightarrow{s}}(E_2/F^{\cyc})$ is cofinitely generated over $\Zp$. Moreover, if this is so, we have the following equality
\[ \la^{\overrightarrow{s}}_{E_1} + \sum_{w\in S'(F^{\cyc})}\corank_{\Zp}\Big(H^1(F^{\cyc}_w,E_1(p))\Big)= \la^{\overrightarrow{s}}_{E_2} + \sum_{w\in S'(F^{\cyc})}\corank_{\Zp}\Big(H^1(F^{\cyc}_w,E_2(p))\Big), \]
where $\la^{\overrightarrow{s}}_{E_i}$ is the $\la$-invariant of the Pontryagin dual of  $\Sel^{\overrightarrow{s}}(E_i/F^{\cyc})$.
\et

\bpf
Clearly, one has $\dim_{\mathbb{F}_p}\Big((E_1(F^{\cyc})[p]\Big) = \dim_{\mathbb{F}_p}\Big((E_2(F^{\cyc})[p]\Big)$ by the congruent assumption. Hence, by Proposition \ref{mod p rank}, it suffices to show that $\Sel^{\overrightarrow{s}}_{non}(E_1[p]/F^{\cyc})\cong \Sel^{\overrightarrow{s}}_{non}(E_2[p]/F^{\cyc})$. Plainly, we have $H^1(G_S(F^{\cyc}), E_1[p])\cong H^1(G_S(F^{\cyc}),E_2[p])$. Therefore, it remains to show that the local terms in the definition of the non-primitive mod-$p$ signed Selmer groups are isomorphic. By $\mathbf{(RED)}$, $E_1$ and $E_2$ have the same reduction type for each prime of $F$ above $p$. Let $\rho$ denote the isomorphism $E_1[p]\cong E_2[p]$. For $w\in S_p^{o}(F^{\cyc})$, write $u$ for the prime of $L'$ below $w$. Denote by $D_{i,u}$ the $\Gal(\overline{L'}_u/L'_u)$-quotient of $E_i(p)$ of $\Zp$-corank one as defined in Subsection \ref{ordinary subsection}.

Now if $w$ is a prime of split multiplicative reduction for $E_1$ and $E_2$, then by Lemma \ref{ordinary}, we have $D_{i,w}(F^{\cyc}_w)\cong\Qp/\Zp$ as $\Gal(\overline{F^{\cyc}_w}/F^{\cyc}_w)$-modules for $i=1,2$. Clearly, in this case, we have $D_{1,w}[p]\cong D_{2,w}[p]$ which in turn implies that $H^1(F^{\cyc}_{w},D_{1,w}[p])\cong H^1(F^{\cyc}_{w}, D_{2,w}[p])$. Now suppose that $w$ is a prime of good ordinary reduction or non-split multiplicative reduction for $E_1$ and $E_2$. Let $u$ be the prime of $L'$ below $w$. Then $u$ is also a prime of good ordinary reduction or non-split multiplicative reduction for $E_1$ and $E_2$.
Since $\rho$ respects the $\Gal(\overline{F}_u/F_u)$-action, $\rho(C_{1,u}[p])$ is a $\Gal(\overline{L'}_u/L'_u)$-submodule of $E_2[p]$ which is isomorphic to $\mu_p$ as $\Gal(\overline{L'}_u/L'^{ur}_u)$-modules. In view of $\mathbf{(RU)}$, we may apply Lemma \ref{D[p]} to conclude that $\rho(C_{1,v}[p]) = C_{2,w}[p]$. Hence $\rho$ induces an isomorphism $D_{1,u}[p]\cong D_{2,u}[p]$ of $\Gal(\overline{L'}_u/L'^{ur}_u)$-modules and hence of $\Gal(\overline{F^{\cyc}_w}/F^{\cyc,ur}_w)$-modules. It then follows that $H^1(F^{\cyc,ur}_{w},D_{1,w}[p])\cong H^1(F^{\cyc,ur}_{w}, D_{2,w}[p])$.

Now suppose that $w\in S_p^{ss}$. As before, write $u$ for a prime of $L'$ below $w$. It follows from hypothesis $\mathbf{(SS)}$ that $L'_u =\Qp$. From the discussion in \cite[Proposition 2.8]{Kim09}, we see that $\rho$ induces an isomorphism $\hat{E}_{1,u}[p]\cong \hat{E}_{2,u}[p]$ of the form $x\mapsto \lambda (ax)$ for some fixed $a\in\mathbb{F}_p$ and some fixed isomorphism $\la$ of formal groups
$\hat{E}_{1,u}\cong \hat{E}_{2,u}$ over $\Zp$. Such an isomorphism in turn induces an isomorphism $H^1(F^{\cyc}_{w},E_1[p])\cong H^1(F^{\cyc}_{w},E_{2}[p])$, which upon restricted to $E^{\pm}_1(F^{\cyc}_w)/p$, gives an isomorphism  $E^{\pm}_1(F^{\cyc}_w)/p\cong E^{\pm}_2(F^{\cyc}_w)/p$. Putting these together, we obtain \[ \frac{H^1(F^{\cyc}_{w},E_1[p])}{E^{\pm}_1(F^{\cyc}_{w})/p}\cong \frac{H^1(F^{\cyc}_{w},E_2[p])}{E^{\pm}_2(F^{\cyc}_{w})/p}.\]
The proof of the theorem is now completed.
\epf

\subsection{$p$-parity conjecture for congruent elliptic curves}

Before proving our next theorem, we first return to the setting of Section \ref{Selmer}, where we only work with one elliptic curve $E$. The following proposition will be required for later discussion. Note that the following proposition does not require $\Sel^{\overrightarrow{s}}(E/F^{\cyc})$ to be cofinitely generated over $\Zp$.

\bp \label{selmer parity} Suppose that $\mathbf{(NA)}$ and $\mathbf{(SS)}$ are valid.
Suppose that $\Sel^{\overrightarrow{s}}(E/F^{\cyc})$ is cotorsion over $\Zp\ps{\Ga}$. In the event that $\overrightarrow{s}\neq \overrightarrow{-}$, assume further that $\mathbf{(S+)}$ holds. Then one has
\[ \corank_{\Zp}\big(\Sel(E/F)\big) = \la^{\overrightarrow{s}}\quad(mod~2),\]
where $\la^{\overrightarrow{s}}$ is the Iwasawa $\la$-invariant of $\Sel^{\overrightarrow{s}}(E/F^{\cyc})$.
\ep

\bpf
When $E$ has good ordinary reduction at all primes above $p$, this was proved in \cite[Proposition 3.10]{G99}. In the case that $E$ has supersingular reduction, this was proved in \cite[Proposition 4.1]{Hat}. The same argument there essentially carries over which we shall sketch briefly. The key to proving the proposition is to make use of the strict signed Selmer groups. By Proposition \ref{strict selmer=selmer} and discussion before the said proposition, the strict signed Selmer group coincides with the usual Selmer group on the level of $F$ and the signed Selmer group on the level of $F^{\cyc}$. Now it is easy to see that the kernel of the natural map
\[ \Sel^{\overrightarrow{s},str}(E/F_n)\stackrel{t_{n}}{\lra} \Sel^{\overrightarrow{s},str}(E/F^{\cyc})^{\Ga_n} \]
is contained in $H^1(\Ga_n, E(F^{\cyc})(p))$. Since the latter is finite and bounded independently of $n$, so is the $\ker t_n$. The point here is to note that in proving this boundness property, we do not require the full strength of a control theorem as in \cite{Maz} or \cite{Kob} which we do not have at our disposal in the context we are working in. The remainder of the proof now proceeds as in \cite[Proposition 3.10]{G99} and \cite[Proposition 4.1]{Hat}.
\epf

\br
In \cite[Theorem 4.6]{AAS}, the authors of that said paper had a result analogous to (but less general than) Proposition \ref{selmer parity}. We like to mention that although their result is correct, it would seem that they did not realize that one has to give a proof via the strict signed Selmer groups as what we have done here (also see \cite[Remark 2.3]{Hat} and \cite[Subsection 4.2]{Kim07}). The reason of using such strict signed Selmer groups is because their local conditions at supersingular primes have a self-dual property in the sense of Proposition \ref{Kim duality} which is required in order to be able to apply Flach's result \cite{Flach}. Finally, it is well-known that the local conditions as defined in the strict signed Selmer groups for other primes have the desired self-dual properties (see \cite[Theorem 1]{Guo}).
\er

The following is an immediate corollary of the previous proposition which generalizes a previous observation of Hatley \cite[Corollary 4.2]{Hat}.

\bc
Suppose that $\mathbf{(NA)}$, $\mathbf{(SS)}$ and $\mathbf{(S+)}$ hold, and that $\Sel^{\overrightarrow{s}}(E/F^{\cyc})$ is torsion over $\Zp\ps{\Ga}$ for every $\overrightarrow{s}$. Then the parity of $\la^{\overrightarrow{s}}$ is independent of $\overrightarrow{s}$.
\ec

\br
We like to take this opportunity to make a remark (which although is not required for subsequent discussion). To the best knowledge of the authors, the signed Selmer groups seem to be treated as independent entities in literature. It is a natural and interesting question to ask whether the signed Selmer groups are related to one another in any way. The above corollary and the result of Hatley \cite[Corollary 4.2]{Hat} seem to suggest (perhaps mildly) that this is so. We also mention another (mild) result in this spirit. In a recent work of the authors \cite[Corollary 2,8]{AL}, we have proved that if the elliptic curve $E$ has good reduction at all primes above $p$ and $\mathbf{(S+)}$ holds, then as long as one of the signed Selmer group vanishes, so will the other (also see \cite[Corollary 5.5]{LLAk} for a generalization of this). It would perhaps be of interest to find out further (possible) relation between the signed Selmer groups.
\er

We come back to the congruent elliptic curves situation.
Recall that the global root number of an elliptic curve $E$ is defined by $w(E/F) = \prod_v w(E/F_v)$, where $v$ runs through all primes of $F$ (see \cite{Do}). The $p$-parity conjecture then predicts that $w(E/F) = (-1)^{s_p(E)}$, where $s_p(E)=\corank_{\Zp}\big(\Sel(E/F)\big)$. We are now in position to prove the following theorem.

\bt \label{p-parity theorem}
Let $E_1$ and $E_2$ be two elliptic curves defined over $L'$ with $E_1[p]\cong E_2[p]$ as $\Gal(\overline{L'}/L')$-modules. Let $L$ be a finite extension of $L'$, and write $F$ for $L$ or $L(\mu_p)$. Suppose that $\mathbf{(NA)}$, $\mathbf{(SS)}$, $\mathbf{(RED)}$ and $\mathbf{(RU)}$ hold for $E_1$ and $E_2$. Suppose further that at least one of the following statements is also valid.
\begin{enumerate}
\item[$(a)$] $\Sel^{\overrightarrow{-}}(E_1/F^{\cyc})$ is cofinitely generated over $\Zp$.
\item[$(b)$] $\Sel^{\overrightarrow{s}}(E_1/F^{\cyc})$ is cofinitely generated over $\Zp$ for some $\overrightarrow{s}\neq \overrightarrow{-}$ and $(\mathbf{S+})$ is valid.
    \end{enumerate}
 Then one has
\[ \frac{w(E_1/F)}{w(E_2/F)} = (-1)^{s_p(E_1) - s_p(E_2)}. \]
In particular, the $p$-parity conjecture holds for $E_1$ over $F$ if and only if it holds for $E_2$ over $F$.
\et

\bpf
Let $v\in S$ be a prime outside $p$. Since $p$ is odd, there are odd number of primes of $S(F^{\cyc})$ above $v$. As $H^1(F^{\cyc}_w, E(p))$ has the common $\Zp$-corank for every $w|v$, it follows that
\[ \corank_{\Zp}\big( H^1(F^{\cyc}_{w_v}, E(p))\big) = \sum_{w|v}\corank_{\Zp}\big( H^1(F^{\cyc}_w, E(p))\big) \quad(\mathrm{mod}~2)\] for some fixed prime $w_v$ of $F^{\cyc}$ above $v$.
Combining this with Theorem \ref{congruent theorem} and Proposition \ref{selmer parity}, we have
\[ s_p(E_1)- s_p(E_2) =\sum_{v\in S'(F)}\big(\si_{E_1, w_v} - \si_{E_2, w_v}\big)\quad (\mathrm{mod}~ 2), \]
where we write $\si_{E_i, w_v}=\corank_{\Zp}\big( H^1(F^{\cyc}_{w_v}, E_i(p))\big)$.

Recall that $w(E/F_v) = -1$ when $v$ is archimedean or is a split multiplicative reduction prime of $E$, and $w(E/F_v) =1$ when $E/F_v$ has good or non-split multiplicative reduction (see \cite[Section 3.4]{Do}). In view of this and that $E_1$ and $E_2$ have the same reduction type for each prime above $p$, we then have
\[ \frac{w(E_1/F)}{w(E_2/F)} = \prod_{v\in S'(F)}\frac{w(E_1/F_v)}{w(E_2/F_v)},\]
 and the latter is equal to
 \[ \prod_{v\in S'(F)}(-1)^{\si_{E_1, w_v} - \si_{E_2, w_v}} \]
  by Proposition \ref{local parity compare}. Combining this with the congruence obtained in the first paragraph of the proof, we obtain the required equality of the theorem. The final assertion of the theorem is then immediate from this.
\epf

\end{document}